\documentclass[11pt]{amsart}

 \usepackage{amsmath,amssymb}
 \usepackage{color}
 \usepackage{mathrsfs}
 \usepackage{graphicx}
 \usepackage{fancybox}
\usepackage{array}

 \numberwithin{equation}{section}
 \theoremstyle{plain}

  {\theoremstyle{definition}}
  
 \newcommand{\reff}[1]{(\ref{#1})}
 
\def\kv{ \mbox{\boldmath$ k$}}
\def\kvi{{\mbox{\footnotesize \boldmath$ k \in \M^j$}\atop \mbox
{\footnotesize \boldmath$|k|$=\it n}}}

 \def\R{ \mathbb R}

 \def\N{{ \mathbb N}}

 \def\M {\mbox{\boldmath$ M$}}

 \def\sp {\quad}

 \def\dis{\displaystyle}

\begin{document}
 
 \title [Derivatives]{The Fa\`a di Bruno formula revisited}


 \author{Raymond Mortini}
  \address{
   \small Universit\'{e} de Lorraine\\
\small D\'{e}partement de Math\'{e}matiques et  
Institut \'Elie Cartan de Lorraine,  UMR 7502\\
\small Ile du Saulcy\\
 \small F-57045 Metz, France} 
 \email{mortini@univ-metz.fr}

 
 \begin{abstract}

{We  present an intuitive approach  to (a variant of) the Fa\`a di Bruno formula which  shows how this formula may have been (re)discovered many times in history.
 Our representation for the $n$-th derivative of the composition $f \circ g$  of two smooth functions
 $f$ and $g$ on $\R$ uses a simpler summation order so that the mysterious condition $b_1+2b_2+\cdots+ nb_n=n$ in Fa\`a di Bruno's formula does not appear. }
 
  \end{abstract}
  
   \maketitle

 
 \section{How one may (re)discover oneself this formula}
 
   Let $f$ and $g$ be two functions on $\R$ whose $n$-th derivatives exist.
   The first 4 derivatives of $f\circ g$ are easy to calculate:   
   
  { \tiny
  $$\begin{matrix}
   (f\circ g)' &= &\hspace{-7mm}(f ' \circ g) g' &&&&&&& \\
     (f\circ g)^{(2)}&=& (f^{(2)}\circ g) (g')^2 &+&\hspace{-13mm}(f '\circ g) g^{(2)} &&&&& \\
   $$ (f\circ g)^{(3)}&= &(f^{(3)}\circ g) (g')^3& +&\hspace{-5mm}3 (f ^{(2)}\circ g) g^{(2)} g'& +& \hspace{-8mm}(f '\circ g) g^{(3)}&&&\\
   (f\circ g)^{(4)}&=& (f^{(4)}\circ g) (g')^4& +& 6 (f^{(3)}\circ g)g^{(2)} (g')^2&+& 
   4 (f^{(2)}\circ g) g^{(3)}g' &+&3 (f^{(2)}\circ g) (g^{(2)})^2&+\;\; ( f '\circ g) g^{(4)}
   \end{matrix}$$
   }

   The expressions get  rapidly longer and more complicated. For  example we have 42 summands for  $n=10$. 
   Note that the number of terms equals the partition number $p(n)$,
   that is the number of ways (without order) of writing the integer $n$ as a sum of strictly positive integers; by the Hardy-Ramanujan formula we have 
   $p(n)\sim \frac{1}{4n\sqrt 3}e^{\pi \sqrt{2n/3}}$.
      
  Let  $$\M^j=\{\mbox{\boldmath$ k$}=(k_1,\dots, k_j)\in (\N^*)^j, \;k_1\geq k_2\geq \cdots\geq k_j\geq 1\}$$
 be the set  of ordered multi-indeces in $\N^*=\{1,2,\dots\}$.
   If $g$ is a function defined on $\R$, and if $g^{(n)}$ is the $n$-th derivative of $g$,
then we denote by $g^{(\mbox{\boldmath$ k$})}$
the function $\prod_{i=1}^j g^{(k_i)}$, where $\mbox{\boldmath$ k$}=(k_1,\dots, k_j)\in \M^j$.
Also,  $g^{(0)}$ is, by convention, equal to the function $g$ itself.

   The classical  Fa\`a di Bruno formula from ca. 1850 gives an explicit formula for $(f\circ g)^{(n)}$:

   $$(f\circ g)^{(n)}(x)=\sum \frac{n!}{b_1!b_2!\cdots b_n!} f^{(j)}(g(x)) \prod_{i=1}^n
   \left(\frac{g^{(i)}(x)}{i!} \right)^{b_i},$$
   where the sum is taken over all different solutions in nonnegative integers $b_1,b_2,\dots, b_n$
   of  $$\mbox{$b_1+2b_2+\cdots +n b_n=n$  and $j:=b_1+\cdots + b_n$.}$$
   A nice historical survey on this appeared in \cite{jp}. See also \cite{cr}. 
   
   Without being aware of that formula,  I developed around 1976-1980  the following  formula:
 
  \begin{equation}\label{mo}
   (f\circ g)^{(n)}(x)= \sum_{j=1}^n f^{(j)}(g(x))\biggl(
 \sum_{\kvi}C_{\mbox{\boldmath$ k$}}^n\, 
g^{(\mathbf  k)}(x)\biggr),
\end{equation}

 where
 $$C_{\mbox{\boldmath$ k$}}^n= \frac{\dis{n\choose \mbox{\boldmath$ k$}} }
 {\dis\prod_i N(\mbox{\boldmath$ k$}, i)!}.$$
Here ${n\choose \mbox{\boldmath$ k$}}$ is  the multinomial coefficient   defined by
 $\dis {n \choose \mbox{\boldmath$ k$}}= \frac{n!}{k_1! k_2!\dots k_j!}$, where
 $|\kv|:=k_1+\dots +k_j=n$, 
and  $N(\mbox{\boldmath$ k$}, i)$ is the number of
times the integer $i$ appears in the  $j$-tuple $\mbox{\boldmath$ k$}$ ($i\in \N^*$ and $\mbox{\boldmath$ k$}\in (\N^*)^j$).

    For example, the coefficient $C^6_{(4,1,1)}$ of   the term
   $ g^{(4)}(g')^2$ when looking at the 6-th  order derivative of $f\circ g$ is
   $$C^6_{(4,1,1)}= \frac{1}{2!}\;\;\frac{  6!}{ 4! \cdot 1! \cdot 1! }= 15$$
   and the coefficient $C^{10}_{(4,2,1,1,1,1)}$ of the  term $g^{(4)} g'' (g')^4$
    in the 10-th derivative is
    $$C^{10}_{(4,2,1,1,1,1)}=\frac{1}{4!}\;\; \frac{10!}{4!\cdot 2!\cdot 1! \cdot 1! 
    \cdot 1! \cdot 1!}   = 3150.$$

  The difference between  our formula and the Fa\`a di Buno formula is that we use a simpler
  summation order and  do not consider exponents of the form $b_j=0$. In particular,  
  we do not need summation over those 
  $(b_1,\dots, b_n)$ satisfying (the difficult to grasp) condition
  $\sum_{i=1}^n i b_i=n$. That these two formulas are 
   equivalent though, immediately follows from a direct comparison of the coefficients.
   Indeed, for fixed $j$ and  
   $\mbox{\boldmath$ k$}=(k_1,\dots,k_j)\in \M^j$, $|\mbox{\boldmath$ k$}|=n$,
    we have:

  $$C_{\mbox{\boldmath$ k$}}^n=\left(\frac {n!}{k_1! \cdots k_j!}\right) \; \frac{1}
  {N(\mbox{\boldmath$ k$}, 1)!}\cdots 
  \frac{1}{N(\mbox{\boldmath$ k$}, n)!}= \biggl(n! \; 
  \frac{1}{\underbrace{i_1!\cdots i_1!}_{b_{i_1} {\rm times}}}\;\;\cdots \;\;
   \frac{1}{\underbrace{i_\ell!\cdots i_\ell!}_{b_{i_\ell} {\rm times}}} \biggr)
   \;\;\frac{1}{b_{i_1}!\cdots b_{i_\ell}!}
  $$
  $$= \frac{n!}{b_1!\cdots b_n!}   \prod_{i=1}^n\left(\frac{1}{i!} \right)^{b_i}$$
   where the $b_{i_m}$ are those exponents that are  different from zero and      where
   $\mbox{\boldmath$ k$}$ has been represented in the canonical form
  $ \mbox{\boldmath$ k$}=(\underbrace{i_1,\dots, i_1}_{{b_{i_1} {\rm times}}}, \cdots,
  \underbrace{i_\ell,\dots, i_\ell}_{{b_{i_\ell} {\rm times}}})$ in decreasing order.
  Note that $$\sum_{i=1}^n ib_i=\sum_{s=1}^\ell i_s b_{i_s} =\sum_{p=1}^j k_p=n$$ and that  
   $$\sum_{i=1}^n{b_i}=b_{i_1}+\cdots +b_{i_\ell}=j.$$
   
   \bigskip

 Next I would like to  present  the (intuitive) steps that led me to the discovery 
 of the formula \reff{mo} above,
   at pre-PC times;  the first (non-programmable) slide rule calculator SR50 had just  appeared.

   1)  I calculated explicitely the derivatives $(f\circ g)^{(n)}$ up to the order 10 
   and wrote them down in a careful chosen order (see figure 1);

2)    An immediate guess is that 

   $$(f\circ g)^{(n)}(x)= \sum_{j=1}^n f^{(j)}(g(x))\biggl(
 \sum_{\kvi}c_{\mbox{\boldmath$ k$}}^n\, 
g^{(\mathbf  k)}(x)\biggr),
$$
for some coefficients  $c_{\mbox{\boldmath$ k$}}^n$ to be determined.

3) Next I gave  an inductive proof that this representation is correct; that needs 
the main step of the construction: where does the factor $g^{(k_1)}\cdots g^{(k_j)}$
with $k_1+\cdots+ k_j=n+1$ comes from? So let us look at the ordered $j$-tuple
$(k_1,\dots, k_j)$. This  tuple is generated, through anti-differentiation, by the $j$ $j$-tuples 
$$(k_1-1,k_2,\dots, k_j),\; (k_1,k_2-1,k_3,\dots,k_j),
\cdots\cdots, (k_1,\dots,k_j-1);$$ 
that is, $$\left(g^{(k_1-1)} g^{(k_2)}\cdots g^{(k_j)}\right)'=
g^{(k_1)}g^{(k_2)}\cdots g^{(k_j)}+\cdots,$$
$$\left(g^{(k_1)} g^{(k_2-1)}\cdots g^{(k_j)}\right)'=
g^{(k_1)}g^{(k_2)}\cdots g^{(k_j)}+\cdots,$$
 etc. 

The main difficulty being that the components $k_j$ are not pairwise distinct.
So their "multiplicities" had to be taken into account.    
This lead to the guess that one may have
   $$c_1^1=1$$
$$ c^{n+1}_{\kv}=\sum_{i=1}^j  \frac{N(\kv -\mathbf e^j_i, k_i-1)}{\dis N(\kv, k_i) }\;
c^n_{\kv -\dis \mathbf e^j_i}$$
   where for $i=1,\dots, j$,  $\mathbf e_i^j=(0,\dots,0,\underbrace{1}_{i-th},0,\dots , 0)$, 
  $\mathbf e_i^j\in \N^j$, $\kv\in \M^j,  |\kv|=n+1$ and $1\leq j\leq n+1$.
  Note that the $j$-tuple $\kv -\mathbf e^j_i$ is not necessarily represented in the canonical form
with decreasing coordinates. Also, if the $i$-th coordinate of $\kv$ is one, then 
the $i$-th coordinate of $\kv -\mathbf e^j_i$ is $0$ and we identify $\kv -\mathbf e^j_i$
with the associated $(j-1)$-tuple. 

 For example in the case $\kv=(3,1,1,1)$
we have
$$ 20=C^6_{(3,1,1,1)}= \frac{1}{1}C^5_{(2,1,1,1)} + \frac{1}{3} C^5_{(3,0,1,1)}+
\frac{1}{3} C^5_{(3,1,0,1)}+\frac{1}{3} C^5_{(3,1,1,0)}, $$
where $ (3, 0,1,1), (3,1,0,1)$ and $(3,1,1,0)$ are identified with  $(3,1,1)$.

With these recursion formula the inductive proof went through. 

4)  Next one has to guess the explicite value of $c_{\mbox{\boldmath$ k$}}^n$, $|\kv|=n$.
Now there are ${n\choose \footnotesize \kv}= \frac{n!}{k_1!\dots k_j!}$ ways to
 choose $k_1$ objects out of $n$, then $k_2$ objects of the remaining ones,  and so on. 
 Due to the multiplicty, one has again to divide by $N(\kv, i)!$.
 
 This gives the guess that $c_{\kv}^n= \frac{\dis{n\choose \mbox{\boldmath$ k$}} }
 {\dis\prod_i N(\mbox{\boldmath$ k$}, i)!}.$
 
5)  These coefficients  $c_{\kv}^n$ actually satisfy the recursion relation above.
Since $c_{1}^1 =1$, and the fact that the recursion relation determines uniquely 
 the next coefficients, we are done: $C_{\kv}^n=c_{\kv}^n$.
\bigskip

\tiny
\renewcommand{\arraystretch}{2}
\setlength{\tabcolsep}{0.1cm}
\hspace{-1,4cm} \begin{tabular} {rcrcccccccccccccc}
 \cline{1-2}
 \multicolumn{1}{|r}{$\mathbf 1$}&\multicolumn{1}{c|}{$f'\,g'$}&&&&&&&&&&&&&&&\\
 \cline{1-4}
  \multicolumn{1}{|r}{$\mathbf 1$}&$f''\,{g'}^2$ &${\mathbf 1}$&\multicolumn{1}{c|}{$f{'}\,{g''}$}&&&&&&&&&&&&&\\
  \cline{1-6}
  \multicolumn{1}{|r}{$\mathbf 1$}&$f^{(3)}\,{g'}^3$&${\mathbf 3}$&$f{''}\,{g''}g'$&$\mathbf {1}$&\multicolumn{1}{c|}{$f{'}\,{g}^{(3)}$}&&&&&&&&&&&\\
  \cline{1-10}
 \multicolumn{1}{|r}{$\mathbf 1$}&$f^{(4)}\,{g'}^4$&${\mathbf 6}$&$f^{(3)}\,{g''}{g'}^2$&$\mathbf {4}$&$f''\,{g}^{(3)}{g'}$&$\mathbf {3}$&$f{''}\,{g''}^2$&$\mathbf {1}$&
\multicolumn{1}{c|}{$f{'}\,{g}^{(4)}$}&&&&&&&\\
\cline{1-16}
  \multicolumn{1}{|r}{$\mathbf 1$} &$f^{(5)}\,{g'}^5$&$\mathbf {10}$&$f^{(4)}\,{g''}{g'}^3$&$\mathbf {10}$&$f^{(3)}\,{g}^{(3)}{g'}^2$&$\mathbf {15}$&$f^{(3)}\,{g''}^{(2)}{g'}$&$\mathbf {5}$&$f{''}\,{g}^{(4)}{g'}$&$\mathbf {10}$&$f{''}\,{g}^{(3)}{g''}$&&&$\mathbf {1}$&\multicolumn{1}{c|}{$f'\,{g}^{(5)}$}&\\
  \cline{1-17}
  \multicolumn{1}{|r}{$\mathbf 1$}&$f^{(6)}\,{g'}^6$&$\mathbf {15}$&$f^{(5)}\,{g''}{g'}^4$&$\mathbf {20}$&$f^{(4)}\,{g}^{(3)}{g'}^3$&$\mathbf {45}$&$f^{(4)}\,{g''}^{(2)}{g'}^2$&$\mathbf {15}$&$f^{(3)}\,{g}^{(4)}{g'}^2$&$\mathbf {60}$&$f^{(3)}\,{g}^{(3)}{g''}{g'}$&$\mathbf {15}$&$\hspace{-6pt}f^{(3)}\,{g''}^3$&$\mathbf {6}$&$f''\,{g}^{(5)}{g'}$&\multicolumn{1}{c|}{$\cdots$}\\
  \cline{1-17}
 \multicolumn{1}{|r}{$\mathbf 1$}&$f^{(7)}\,{g'}^7$&$\mathbf {21}$&$f^{(6)}\,{g''}{g'}^5$&$\mathbf {35}$&$f^{(5)}\,{g}^{(3)}{g'}^4$&$\mathbf {105}$&$f^{(5)}\,{g''}^{(2)}{g'}^3$&$\mathbf {35}$&$f^{(4)}\,{g}^{(4)}{g'}^3$&$\mathbf {210}$&$f^{(4)}\,{g}^{(3)}{g''}{g'}^2$&$\mathbf {105}$&$f^{(4)}\,{g''}^3{g'}$&$\mathbf {21}$&$f^{(3)}\,{g}^{(5)}{g'}^2$&\multicolumn{1}{c|}{$\cdots$}\\
   \cline{1-17}
 \multicolumn{1}{|r}{$\mathbf 1$}&$f^{(8)}\,{g'}^8$&$\mathbf {28}$&$f^{(7)}\,{g''}{g'}^6$&$\mathbf {56}$&$f^{(6)}\,{g}^{(3)}{g'}^5$&$\mathbf {210}$&$f^{(6)}\,{g''}^{(2)}{g'}^4$&$\mathbf {70}$&$f^{(5)}\,{g}^{(4)}{g'}^4$&$\mathbf {560}$&$f^{(5)}\,{g}^{(3)}{g''}{g'}^3$&$\mathbf {420}$&$f^{(5)}\,{g''}^3{g'}^2$&$\mathbf {56}$&$f^{(4)}\,{g}^{(5)}{g'}^3$&\multicolumn{1}{c|}{$\cdots$}\\
  \cline{1-17}
 \multicolumn{1}{|r}{$\mathbf 1$}&$f^{(9)}\,{g'}^9$&$\mathbf {36}$&$f^{(8)}\,{g''}{g'}^7$&$\mathbf {84}$&$f^{(7)}\,{g}^{(3)}{g'}^6$&$\mathbf {378}$&$f^{(7)}\,{g''}^{(2)}{g'}^5$&$\mathbf {126}$&$f^{(6)}\,{g}^{(4)}{g'}^5$&$\mathbf {1260}$&$f^{(6)}\,{g}^{(3)}{g''}{g'}^4$&$\mathbf {1260}$&$f^{(6)}\,{g''}^3{g'}^3$&$\mathbf {126}$&$f^{(5)}\,{g}^{(5)}{g'}^4$&\multicolumn{1}{c|}{$\cdots$}\\
   \cline{1-17}
 \multicolumn{1}{|r}{$\mathbf 1$}&$f^{(10)}\,{g'}^{10}$&$\mathbf {45}$&$f^{(9)}\,{g''}{g'}^8$&$\mathbf {120}$&$f^{(8)}\,{g}^{(3)}{g'}^7$&$\mathbf {630}$&$f^{(8)}\,{g''}^{(2)}{g'}^6$&$\mathbf {210}$&$f^{(7)}\,{g}^{(4)}{g'}^6$&$\mathbf {2520}$&$f^{(7)}\,{g}^{(3)}{g''}{g'}^5$&$\mathbf {3150}$&$f^{(7)}\,{g''}^3{g'}^4$&$\mathbf {252}$&$f^{(6)}\,{g}^{(5)}{g'}^5$&\multicolumn{1}{c|}{$\cdots$}\\
\hline
 \end{tabular}

\normalsize

\bigskip

\section{Further formulas  and questions}

Applying  formula \reff{mo} for the function $f(x)=\log x, x>0$ and $g(x)=e^x$
gives
$$\sum_{j=1}^n (-1)^{j-1}(j-1)!    \hspace{-10pt}\sum_{\kvi}C_{\mbox{\boldmath$ k$}}^n=0, \sp n\geq 2\sp ;$$

whereas for $f(x)=x^n$ and $g(x)=e^x$ one obtains
$$\sum_{j=1}^n {n\choose j} j!   \hspace{-5pt}\sum_{\kvi}C_{\mbox{\boldmath$ k$}}^n=n^n.$$
In particular, $L:=\dis \sum_{\footnotesize\kv: |\kv| =n}C^n_{\kv}\leq n^n$.
Is there an explicit expression for $L$? If one uses $f(x)=g(x)=e^x$, then
$$\left(e^{e^x}\right)^{(n)}\left|_{x=0}\right.=eL.$$

One may also ask the  following questions:\medskip

(1) What is $\dis \sum_{\kvi} C^n_{\kv}$\;\; $(1\leq j\leq n)$?

(2) What is $\max \{C^n_{\kv}: |\kv|=n\}$?\medskip

(3) Is there a formula for the number of partitions of $n$ with fixed length $j$?

\medskip

In our scheme (figure 1), one can give easy formulas for the coefficients in each column.
In fact, each element in  a fixed column is a multiple of the first coefficient. More precisely,
if $k_1\geq k_2\geq k_j>1$, then we have:
$$ C^{n+\ell}_{(k_1,\dots,k_j, \underbrace{1,\dots, 1}_{\ell-\text{times}})}=
C^n_{(k_1,\dots,k_j)} {n+\ell\choose \ell},\;\;
 \ell\in \N.$$
 
 We observe that several columns coincide; for example 
 $C^6_{(2,2,2)}=C^6_{(4,2)}$ and so the elements of the associated columns are the same.
 
 There are actually infinitely many pairs of colums that coincide (just use
 that $ C^{6+i}_{(2,2,2,i)}= C^{6+i}_{(4,2,i)}$ for every $i\geq 5$.)
 
 Are there triples (or higher number)  of columns  that coincide?

\bigskip

{\bf Acknowledgements} I thank J\'er\^ome No\"el for his help in finding the exact 
{\TeX}-commands
for creating the  tabular form above. I also thank Claude Jung for valuable discussions
around 1980.

 \end{document}